# ON THE HYPERBOLICITY OF THE FEIGENBAUM FIXED POINT

DANIEL SMANIA

ABSTRACT. We show the hyperbolicity of the Feigenbaum fixed point using the inflexibility of the Feigenbaum tower, the Manẽ-Sad-Sullivan $\lambda$- Lemma and the existence of parabolic domains (petals) for semi-attractive fixed points.

## 1. INTRODUCTION AND STATEMENT OF RESULTS

Let $g\colon U_0 \to V_0$ be a quadratic-like map. This means that $g$ is a ramified holomorphic covering map of degree two, where $U$ and $V$ are simply connected domains, $U \Subset V$. We also assume that the filled-in Julia set of $g$, $K(g) := \cap_n g^{-n} V$, is connected. We say that $g$ is renormalizable with period two if there exist simply connected subdomains $U_1$, $V_1$ so that $g^2 \colon U_1 \to V_1$ is also a quadratic-like map with connected filled-in Julia set.

Two quadratic like maps $h\colon U_h \to V_h$ and $g\colon U_g \to V_g$, both with connected filled-in Julia set, defines the same quadratic-like germ if $K(h)$ coincides with $K(g)$ and $h$ coincides with $g$ in a neighborhood of $K(g)$. If $g$ is renormalizable, then the renormalization of the germ defined by $g$ is the unique quadratic-like germ defined by the normalization of any possible induced map $g^2 \colon U_1 \to V_1$ which are quadratic-like maps with connected filled-in Julia set (normalize the germ using an affine conjugacy, setting the critical point at zero and the unique fixed point in $K(g_1)$ which does not cut $K(g_1)$ in two parts, the so-called $\beta$ fixed point of $g_1$, to 1). The operator $\mathcal{R}$ is called the Feigenbaum renormalization operator. In the setting of quadratic-like germs which have real values in the real line, there exists an unique fixed point to the Feigenbaum renormalization operator (proved by D. Sullivan: see also [McM96]), denoted $f^\star$ (it is a open question if this is the unique fixed point in the set of all quadratic-like germs).

It is a consequence of the so-called a priori bounds [Lyu99] that we can choose a simply connected domain $U$, $K(f) \subset U$, so that, if $\mathcal{B}(U)$ denotes the Banach space of the complex analytic functions $g$, $Dg(0) = 0$, with a continuous extension to $\overline{U}$, provided with the sup norm, and $\mathcal{B}_{nor}(U)$ denotes the affine subspace of the functions $g$ so that $g(1) = 1$, then the fixed point $f^\star$ has a complex analytic extension which belongs to $\mathcal{B}_{nor}(U)$ and there exists $N$ so that the operator $\mathcal{R}^N$ can be represented as a compact operator defined in a small neighborhood of $f^\star$ in

2000 *Mathematics Subject Classification.* 37F25 , 37F45, 37E20.

This work was partially supported by CNPq-Brazil grant 200764/01-2 and University of Toronto. I would like to thank IMS-SUNY at Stony Brook, University of Toronto and specially to M. Lyubich by my wonderful stay at Stony Brook and Toronto. I am also in debit to M. Lyubich by the useful comments.

$\mathcal{B}_{nor}(U)$. More precisely, there exists a larger domain $\tilde{U} \supsetneq U$ and a complex analytic operator $\tilde{\mathcal{R}}\colon B_{\mathcal{B}_{nor}(U)}(f^\star, \epsilon) \to \mathcal{B}_{nor}(\tilde{U})$ so that, if $i$ denotes the natural inclusion $i\colon \mathcal{B}(\tilde{U}) \to \mathcal{B}(U)$, then $\mathcal{R}^N = i \circ \tilde{\mathcal{R}}$, where the equality holds in the intersection of the domains of the operators. To simplify the notation, we will assume that $N = 1$ and identify $\mathcal{R}$ with its complex analytic extension in $\mathcal{B}_{nor}(U)$.

Two quadratic-like maps $g_0$ and $g_1$ are in the same hybrid class if there exists a quasiconformal conjugacy $\phi$ between them, in a neighborhood of their filled-in Julia sets, so that $\overline{\partial}\phi \equiv 0$ on $K(g_0)$. Note that quadratic-like maps in the hybrid class of $f^\star$ are infinitely renormalizable.

We will provide a new approach to the following result:

**Theorem 1** (Exponential contraction:[McM96] and [Lyu99]). *There exists $\lambda < 1$ so that, for every quadratic-like map $f$ which is in the hybrid class of $f^\star$, there exist $n_0 = n_0(f)$ and $C = C(f) > 0$ so that $\mathcal{R}^n f \in B_{\mathcal{B}(U)}(f^\star, \epsilon)$, for $n \geq n_0$, and $|\mathcal{R}^{n_0+n} f - f^\star|_{\mathcal{B}(U)} \leq C\lambda^n$.*

A major attractive of this new proof is that it is essentially infinitesimal and has a "dynamical flavor": we will prove that the derivative of the renormalization operator is a contraction in the tangent space of the hybrid class (the contraction of the derivative of the renormalization operator on the hybrid class was proved by Lyubich [Lyu99], but his proof is not infinitesimal). Moreover, the method seems to be so general as the previous ones: it also applies to the classical renormalization horseshoe [Lyu99] and the Fibonacci renormalization operator [Sm02a], for instance.

We will also obtain, as a corollary of McMullen theory of towers [McM96], the local behavior of semi-attractive fixed points [H] and an easy application of the $\lambda$-lemma [MSS] that

**Theorem 2** ([Lan][Lyu99]). *The Feigenbaum fixed point is hyperbolic.*

The reader will observe that we assume the Feigenbaum combinatorics just to simplify the notation: the argument in the proof of Theorem 2 works as well to prove the hyperbolicity of real periodic points of the renormalization horseshoe.

## 2. Preliminaries

2.1. **Parabolic domains for semi-attractive fixed points.** Consider a complex Banach space $B$, and let $F\colon A \subset B \to B$ be a complex analytic operator defined in an open set $A$. Suppose that $p \in A$ is a fixed point for $F$. We say that $p$ is a **semi-attractive** fixed point for $F$ if

- The value 1 is an eigenvalue for $DF_p$.
- There exists a Banach subspace $E^s$, with (complex) codimension one, which is invariant by the action of $DF_p$ and furthermore the spectrum of $DF_p$, restricted to $E^s$, is contained in $\{z\colon |z| \leq r\}$, where $r < 1$.

The following result was proved by M. Hakim [H] for finite-dimensional complex Banach spaces ($\mathbb{C}^n$), but the proof can be carry out as well for a general complex Banach space:

**Proposition 2.1** ([H]). *Consider a compact complex analytic operator $F$, defined in an open set of a complex Banach space $B$. Let $p$ be a semi-attractive fixed point. Then one of the following statements holds:*

(1) **Curve of fixed points:** *There exists a complex analytic curve of fixed points which contains p.*
(2) **Parabolic domains (Petals):** *There exists $k \geq 1$ so that, for every $\epsilon > 0$ there exists a connected open set $U$, whose diameter is smaller than $\epsilon$, which is forward invariant by the action of $F$ and, moreover,*

$$F^n u \to_n p, \text{ for every } u \in U,$$

*where the speed of this convergence is subexponential: for each $u \in U$, there exists $C = C(u)$ so that*

$$\frac{1}{C}\frac{1}{n^{1/k}} \leq |F^n u - p| \leq C\frac{1}{n^{1/k}}.$$

An outline of Hakim's proof can be found in the Appendix.

## 3. Infinitesimal contraction on the horizontal space

Let $f\colon V_1 \to V_2$ be a quadratic-like map with connected Julia set and with an analytical extension to $\mathcal{B}_{nor}(U)$, with $K(f) \subset U$. The horizontal subspace (introduced by Lyubich[Lyu99]) of $f$, denoted $E_f^h$, is the subspace of the vectors $v \in \mathcal{B}(U)$ so that there exists a quasiconformal vector field in the Riemann sphere $\alpha$ satisfying $v = \alpha \circ f - Df \cdot \alpha$ in a neighborhood of $K(f)$, with $\overline{\partial}\alpha \equiv 0$ on $K(f)$ and $\alpha(0) = \alpha(1) = \alpha(\infty) = 0$. We will not use the following information here, but certainly it will clarify the spirit of our methods: in an appropriated setting, the hybrid class is a complex analytic manifold and the horizontal space is the tangent space of the hybrid class at $f$(see [Lyu99]).

**Lemma 3.1** ([ALdM]). *Let $f$ be a quadratic-like map with an extension to $\mathcal{B}_{nor}(U)$ and connected Julia set contained in $U$. Assume that $f$ does not support invariant line fields in its filled-in Julia set. Let $V \Subset U$ be a domain with smooth boundary so that $K(f) \subset V$. Then there exist $C, \epsilon > 0$ so that, if $|f - g|_{\mathcal{B}(U)} \leq \epsilon$ and $g\colon g^{-1}V \to V$ is a quadratic-like map with connected Julia set, then, for every $v \in E_g^h$ there exists a $C|v|_{\mathcal{B}(U)}$-quasiconformal vector field $\alpha$ in $\overline{\mathbb{C}}$ so that $v = \alpha \circ g - Dg \cdot \alpha$ on $V$.*

With the aid of a compactness criterium to quasiconformal vector fields in $\overline{\mathbb{C}}$, we have:

**Corollary 3.2** ([ALdM]). *Assume that $(f_n, v_n) \to_n (f_\infty, v_\infty)$ in $\mathcal{B}_{nor}(U) \times \mathcal{B}(U)$, where $f_i \colon f_i^{-1} V \to V$, $i \in \mathbb{N} \cup \{\infty\}$, are quadratic-like maps with connected filled-in Julia sets $K(f_i) \subset V \Subset U$. Furthermore, assume that $v_n \in E_{f_n}^h$, for $n \in \mathbb{N}$. If $f$ does not support invariant line fields in $K(f)$, then $v_\infty \in E_f^h$. In particular $E_f^h$ is closed.*

If $\mathcal{R}$ is the nth iteration of the Feigenbaum renormalization operator and $f$ is close to $f^\star$ in $\mathcal{B}(U)$, denote by $\beta_f$ the analytic continuation of the $\beta$-fixed point of the small Julia set associated with the nth renormalization of $f^\star$ [McM96]. The following result gives a description of the action of the derivative in a horizontal vector $v = \alpha \circ f - Df \cdot \alpha$ in terms of $\alpha$:

**Proposition 3.3.** *Let $V$ be a neighborhood of $K(f^\star)$. Replacing $\mathcal{R}$ by an iteration of it, if necessary, the following property holds: If $f \in \mathcal{B}_{nor}(U)$ is close enough to* $f^\star$ and $v = \alpha \circ f - Df \cdot \alpha$ on $V$, where $v \in \mathcal{B}(U)$ and $\alpha$ is a quasiconformal vector field in the Riemann sphere, normalized by $\alpha(0) = \alpha(1) = \alpha(\infty) = 0$, then

$$D\mathcal{R}_f \cdot v = r(\alpha) \circ \mathcal{R}f - D(\mathcal{R}f) \cdot r(\alpha), \tag{1}$$

on $U$, where

$$r(\alpha)(z) := \frac{1}{\beta_f}\alpha(\beta_f z) - \frac{1}{\beta_f}\alpha(\beta_f) \cdot z.$$

In particular, if $f$ is renormalizable, then $D\mathcal{R}_f E_f^h \subset E_{\mathcal{R}f}^h$.

This result is consequence of a simply calculation and the complex bounds to $f^\star$. Note that, apart the normalization by a linear vector field, $r(\alpha)$ is just the pullback of the vector field $\alpha$ by a linear map. In particular, if $\alpha$ is a $C$-quasiconformal vector field, then $r(\alpha)$ is also a $C$-quasiconformal vector field: this will be a key point in the proof of the infinitesimal contraction of the renormalization operator in the horizontal subspace (Proposition 3.4).

Let $f^\star\colon V_1 \to V_2$ be a quadratic-like representation of the fixed point. The Feigenbaum tower is the indexed family of quadratic-like maps $f_i^\star\colon \beta_{f^\star}^{-i} V_1 \to \beta_{f^\star}^{-i} V_2$, $i \in \mathbb{N}$, defined by $f_i^\star(z) := \beta_{f^\star}^{-i} f^\star(\beta_{f^\star}^i \cdot z)$.

**Proposition 3.1** ([McM96]). *The Feigenbaum tower does not support invariant line fields: this means that there is not a measurable line field which is invariant by all (or even an infinite number of) maps in the Feigenbaum tower.*

**Proposition 3.2** ([Su] and [McM96]). *Let $f$ be a quadratic-like map which admits a hybrid conjugacy $\phi$ with $f^\star$. Then $\phi_n(z) := \beta_{f^\star}^{-n} \cdot \phi(\beta_{\mathcal{R}^{n-1}f} \cdots \beta_f \cdot z)$ converges to identity uniformly on compact sets in the complex plane. In particular, there exists $n_0 = n_0(f)$ so that $\mathcal{R}^n f \in B_{\mathcal{B}_{nor}(U)}(f^\star, \epsilon)$, for $n > n_0$, and $\mathcal{R}^n f \to_n f^\star$ on $\mathcal{B}_{nor}(U)$.*

Theorem 1 says that this convergence is, indeed, exponentially fast. The following proposition has a straightforward proof:

**Proposition 3.3.** *Let $\mathcal{R} = i \circ \tilde{\mathcal{R}}$, where $\tilde{\mathcal{R}}\colon \mathcal{V} \subset \mathcal{B} \to \tilde{\mathcal{B}}$ is an operator defined in a neighborhood $\mathcal{V}$ of a Banach space $\mathcal{B}$, to another Banach space $\tilde{\mathcal{B}}$, and $i\colon \tilde{\mathcal{B}} \to \mathcal{B}$ is a compact linear transformation. Let $S \subset \mathcal{B} \times \mathcal{B}$ be a set with the following properties:*

(1) **Vector bundle structure:** *If $(f, v_1)$ and $(f, v_2) \in S$, then $(f, \alpha \cdot v_1 + v_2) \in S$, for every $\alpha \in \mathbb{C}$,*
(2) **Semicontinuity:** *If $(f_n, v_n) \to (f, v)$ and $(f_n, v_n) \in S$, then $(f, v) \in S$,*
(3) **Invariance:** *If $(f, v) \in S$ then $(\mathcal{R}f, D\mathcal{R}_f \cdot v) \in S$,*
(4) **Compactness:** *$\{(\tilde{\mathcal{R}}f, D\tilde{\mathcal{R}}_f \cdot v)\colon (f, v) \in S, |v| \leq 1\}$ is a bounded set in $\tilde{\mathcal{B}} \times \tilde{\mathcal{B}}$,*
(5) **Uniform continuity:** *Denote $E_f := \{(f, v)\colon (f, v) \in S\}$. There exists $C > 0$ so that, for every $f$ and $n \geq 0$, $|D\mathcal{R}_f^n|_{E_f} \leq C$,*
(6) *If $(f, v) \in S$ then $|D\mathcal{R}_f^n \cdot v| \to_n 0$,*

*Then there exist $\lambda < 1$ and $N \in \mathbb{N}$ so that $|D\mathcal{R}_f^N|_{E_f} \leq \lambda$, for every $f$ so that $E_f \neq \varnothing$.*

Proposition 3.3 is a generalization of the following fact about compact linear operators $T\colon \mathcal{B} \to \mathcal{B}$: if $T^n v \to 0$, for every $v \in \mathcal{B}$, then the spectral radius of $T$ is strictly smaller than one.

Given $\epsilon$, $K > 0$, and a domain $V \Subset U$ so that $K(f^\star) \subset V$, denote by $\mathcal{A}(\epsilon, K, V)$ the set of maps $f \in \mathcal{B}_{nor}(U)$ so that there exists a $K$- quasiconformal map $\phi$ in the complex plane so that $\phi(V) \subset \overline{U}$ and $\phi \circ f^\star = f \circ \phi$ on $V$; moreover, for $n \geq 0$, we have $|\mathcal{R}^n f - f^\star|_{\mathcal{B}(U)} \leq \epsilon$. Note that $\mathcal{A} := \mathcal{A}(\epsilon, K, V)$ is closed. Furthermore, replacing $\mathcal{R}$ by an iterate, if necessary, we can assume that $\mathcal{A}$ is invariant by the action of $\mathcal{R}$. Selecting $K$ and $\epsilon$ properly, by the topological convergence (Proposition 3.2) and Lemma 2.2 in [Lyu02], for every $f$ in the hybrid class of $f^\star$, there exists $N = N(f)$ so that $\mathcal{R}^N f \in \mathcal{A}$.

**Proposition 3.4** (Infinitesimal contraction: cf. [Lyu99])**.** *There exist $\lambda < 1$ and $N > 0$ so that $|D\mathcal{R}^N_f|_{E^h_f} \leq \lambda$, for every $f \in \mathcal{A}(\epsilon, K, V)$.*

*Proof.* Consider the set $S := \{(f, v) \colon f \in \mathcal{A}, \ v \in E^h_f\}$. It is sufficient to verify the properties in the statement of Proposition 3.3. Since $\mathcal{A}$ is closed, property 2 follows of Corollary 3.2. Since $\mathcal{A}$ is invariant by $\mathcal{R}$, property 3 follows of Proposition 3.3. The compactness property is obvious, if $\epsilon$ is small enough. To prove the uniform continuity property, by Propositions 3.1 and 3.3, we have that, for $(f, v) \in S$ and $n \geq 1$, $D\mathcal{R}^n_f \cdot v = \alpha_n \circ \mathcal{R}^N f - D(\mathcal{R}^N f) \cdot \alpha_n$ on $U$, with

$$\alpha_n(z) := \frac{1}{\beta_{n-1} \ldots \beta_0} \alpha(\beta_{n-1} \ldots \beta_0 z) - \frac{1}{\beta_{n-1} \ldots \beta_0} \alpha(\beta_{n-1} \ldots \beta_0) z,$$

where $\beta_i = \beta_{\mathcal{R}^i f}$ and $\alpha_n$ are $K \cdot |v|_{\mathcal{B}(U)}$-quasiconformal vector fields. Note that $K$ does not depends on $(f, v) \in S$ or $n \geq 1$. By the compactness of $K$- quasiconformal vector fields (recall that $\alpha_n(0) = \alpha_n(1) = \alpha_n(\infty) = 0$), we get $|D\mathcal{R}^n_f|_{E^h_f} \leq C$, for some $C > 0$. To prove assumption 6, note that $\overline{\partial}\alpha_n$ is an invariant Beltrami field to the finite tower

$$\mathcal{R}^n f, \frac{1}{\beta_{n-1}} \mathcal{R}^{n-1} f(\beta_{n-1} z), \ldots, \frac{1}{\beta_{n-1} \cdots \beta_0} f(\beta_{n-1} \cdots \beta_0 z).$$

But, by the topological convergence, these finite towers converges to the Feigenbaum tower. Hence, if a subsequence $\alpha_{n_k}$ converges to a quasiconformal vector field $\alpha_\infty$, then $\overline{\partial}\alpha_\infty$ is an invariant Beltrami field to the Feigenbaum tower (since $\overline{\partial}\alpha_{n_k}$ converges to $\overline{\partial}\alpha_\infty$ in the distributional sense), so, by Proposition 3.1, $\alpha_\infty$ is a conformal vector field in the Riemann sphere. Since $\alpha_\infty$ vanishes at three points, $\alpha_\infty \equiv 0$. Hence $\alpha_n \to 0$ uniformly on compact sets in the complex plane, so we get $\mathcal{R}^n_f \cdot v \to 0$ (Note that $|D(\mathcal{R}^n f)|$ is uniformly bounded, for $n \geq 1$). $\square$

We are going to prove Theorem 1: Let $f$ be a quadratic-like map in the hybrid class of $f^\star$. Then there exists a quasiconformal map $\phi \colon \mathbb{C} \to \mathbb{C}$ which is a conjugacy between them in a neighborhood of their Julia sets. Consider the following Beltrami path $f_t$ between the two maps, induced by $\phi$: if $\phi_t$, $|t| \leq 1$, is the unique normalized quasiconformal map so that $\overline{\partial}/\partial \phi_t = t \cdot \overline{\partial}/\partial \phi$, then $f_t = \phi_t \circ f \circ \phi_t^{-1}$. By the topological convergence, there exists $n_0$ so that $\mathcal{R}^{n_0+n} f_t \in \mathcal{A}$, for $n \geq 0$, $|t| \leq 1$. An easy calculation shows that

$$\left. \frac{d\mathcal{R}^{n_0+n} f_t}{dt} \right|_{t=t_0} \in E^h_{\mathcal{R}^{n_0+n} f_{t_0}},$$

for $|t_0| \leq 1$. The infinitesimal contraction finishes the proof.

*Remark* 1. *The first step in the above proof on Theorem 1, to prove that $\alpha_n \to 0$ (in the proof of Proposition 3.4), must be compared with the proof of Lemma*

*9.12 in* [McM96]. *In C. McMullen argument, additional considerations should be done to arrive in exponential contraction; firstly it is proved that quasiconformal deformations (as the quasiconformal vector field $\alpha$ in the definition of the horizontal vectors) are $C^{1+\beta}$-conformal at the critical point ( Lemma 9.12 in* [McM96] *and the deepness of the critical point have key roles in this proof), and then it is necessary to integrate this result. In M. Lyubich argument* [Lyu99], *firstly it is proved that the hybrid class is a complex analytic manifold and then the topological convergence is converted in exponential contraction via Schwartz's Lemma.*

## 4. Hyperbolicity of the Feigenbaum fixed point

We are going to prove Theorem 2. Firstly we will prove that

(2)  $$\sigma(D\mathcal{R}_{f^\star}^2) \cap \mathbb{S}^1 \subset \{1\}.$$

Indeed, if $D\mathcal{R}_{f^\star} \cdot v = \lambda v$, then the vector $\tilde{v}(z) := \overline{v(\overline{z})}$ is a solution to $D\mathcal{R}_{f^\star} \cdot \tilde{v} = \overline{\lambda}\tilde{v}$. So if $\lambda \in \mathbb{S}^1 \setminus \{-1, 1\}$ then *codim* $E^h > 1$, which is a contradiction ($E^h$ has codimension one [Lyu99]. The same result can be proven in an easy way using the argument explained in section 12 on [Sm02a]). Indeed, we can prove, using the contraction on the horizontal direction and results on [Sm02b], which uses only elementary methods, that $\sigma(D\mathcal{R}_{f^\star}) \cap \mathbb{S}^1 \subset \{1\}$, but the proof is more involving.

Furthermore $\sigma(D\mathcal{R}_{f^\star})$ is not contained in $\mathbb{D}$ (see Lyubich[Lyu99]. We can also use the results in [Sm02b] to prove this claim). So either $f^\star$ is a hyperbolic fixed point (with a onedimensional expanding direction) or it is a semi-attractive fixed point, since by Proposition 3.4 the derivative of the renormalization operator at the fixed point is a contraction on the horizontal space, which has codimension one. Assume that $f^\star$ is semi-attractive and let's arrive in a contradiction. Indeed, by Proposition 2.1, one of the following statements holds:

**Case i.** There exists a connected open set of maps $\mathcal{U} \subset \mathcal{B}_{nor}(U)$, whose diameter can be taken small, which is forward invariant by the action of $\mathcal{R}^2$ and so that each map in $\mathcal{U}$ is attracted at a subexponential speed to the fixed point $f^\star$. Because the maps in $\mathcal{U}$ are very close to $f^\star$ and $\mathcal{U}$ is forward invariant, all the maps in $\mathcal{U}$ are infinitely renormalizable (this argument is easy: see Lemma 5.8 in [Lyu99]). So their filled-in Julia sets have empty interior and their periodic points are repelling, hence there are not bifurcations of periodic points in $\mathcal{U}$. Consider two maps $g, \tilde{g}$ in $\mathcal{U}$ which admit a complex analytic path $g \colon \mathbb{D} \to \mathcal{U}$ between them ($g_0 = g$ and, for some $|\lambda| < 1$, $\tilde{g} = g_\lambda$). Because $\mathbb{D}$ is simply connected and there are not bifurcations of periodic points in $\mathcal{U}$, each periodic point $p \in K(g_0)$ has an unique analytic continuation $h(p, \lambda)$, $\lambda \in \mathbb{D}$: this means that $h(p, 0) = p$ and $h(p, \lambda)$ is a periodic point of $g_\lambda$. So the function $h \colon Per(g_0) \times \mathbb{D} \to \mathbb{C}$ defines a holomorphic motion on $Per(g_0) = \{p \colon \exists n > 0 \text{ s.t. } g_0^n(p) = p\}$ (note that $h(p, \lambda) \neq h(q, \lambda)$, if $p \neq q$, since there are not bifurcations of periodic points). Moreover, provided $\mathcal{U}$ is small enough, we can select a domain $U_1$ with a real analytic Jordan curve boundary so that, for every $\lambda \in \mathbb{D}$, $g_\lambda \colon g_\lambda^{-1} U_1 \to U_1$ is a quadratic-like representation. We can also easily define a holomorphic motion $h \colon U_1 \setminus g_0^{-1} U_0 \times \mathbb{D} \to \mathbb{C}$ so that $h(x, \lambda) \equiv x$, for $x \in \mathbb{C} \setminus U_1$ and $g_\lambda(h(x, \lambda)) = h(g_0(x), \lambda)$, for $x \in \partial g_0^{-1} U_1$. Since $g_\lambda$ have connected filled-in Julia sets, we can extend the holomorphic motion to a holomorphic motion $h \colon \mathbb{C} \setminus K(g_0) \times \mathbb{D} \to \mathbb{C}$ so that $g_\lambda(h(x, \lambda)) = h(g_0(x), \lambda)$, for $x \in g_0^{-1} U_1 \setminus K(g_0)$. So we have defined a holomorphic motion $h$ on the everywhere dense set $Per(g_0) \cup (\mathbb{C} \setminus K(g_0))$ which commutes with the dynamics. By the $\lambda$-lemma [MSS], this holomorphic motion extends to the whole Riemann sphere, so all maps $g_\lambda$ are quasiconformaly conjugated. Since there is a piecewise complex analytic path between any two maps in $\mathcal{U}$, we conclude that all maps in $\mathcal{U}$ are in the same quasiconformal class. Note that the above construction does not give any upper bound for the quasiconformality of the conjugacy: the quasiconformality could be large when the Kobayashi distance between $g$ and $\tilde{g}$ on $\mathcal{U}$ is large.

We claim that, provided $\mathcal{U}$ is small enough, it is possible to choose a quasiconformal conjugacy between any two maps in $\mathcal{U}$ so that the quasiconformality is uniformly bounded *outside* their filled-in Julia sets, using the argument in the proof of Lemma 2.3 in [Lyu02]: in a small neighborhood $\mathcal{V} \subset \mathcal{B}_{nor}(U)$ of $f^\star$, it is possible to find a domain $U_1$ so that $g\colon g^{-1}U_1 \to U_1$ is a quadratic-like restriction of $g$ (but note that the Julia sets of these quadratic-like restrictions are not, in general, connected). This define the holomorphic moving fundamental annulus $U_1 \setminus g^{-1}U_1$. In particular, provided $\mathcal{U}$ is small enough, there exists $B > 0$ so that for every $g_0$ and $g_1$ which belong to $\mathcal{U}$, there exists a $B$-quasiconformal mapping $h$ between $\mathbb{C}\setminus g_0^{-1}U_1$ and $\mathbb{C} \setminus g_1^{-1}U_1$ so that $h \equiv Id$ on $\mathbb{C} \setminus U_1$ and $g_1 \circ h = h \circ g_0$ on $\partial g_0^{-1}U_1$. Since the Julia sets of $g_0$ and $g_1$ are connected, we can extend $h$ to a $B$-quasiconformal map

$$h\colon \mathbb{C} \setminus K(g_0) \to \mathbb{C} \setminus K(g_1)$$

which is a conjugacy on $g_0^{-1}U_1 \setminus K(g_0)$. Once we already know that $g_0$ and $g_1$ are in the same quasiconformal class, $h$ has a quasiconformal extension $h_{g_0,g_1}$ to $\mathbb{C}$ (this follows as in the proof of Lemma 1, in [DH, pg. 302]: if $\tilde{h}$ is a quasiconformal conjugacy between $g_0$ and $g_1$, then $\tilde{h}^{-1} \circ h$ commutes with $g_0$ outside $K(g_0)$, which implies that $\tilde{h}^{-1} \circ h$ extends to a homeomorphism in $\mathbb{C}$ which coincides with $Id$ on $K(g_0)$. By the Rickmann removability theorem (see the statement in [DH]), this map is a quasiconformal homeomorphism, so $h$ is a quasiconformal homeomorphism). This finishes the proof of the claim.

Since all renormalizations of these maps are very close to $f^\star$, they also satisfies the unbranched complex bounds condition (this is consequence of a short lemma in [Lyu97]). In particular there are not invariant line fields supported on their filled-in Julia sets [McM94], and hence the quasiconformality of the conjugacy $h_{g_0,g_1}\colon \mathbb{C} \to \mathbb{C}$ is uniformly bounded on the whole complex plane by $B$. But $f^\star$ is a boundary point of $\mathcal{U}$, so the compactness of $B$-quasiconformal maps (note that the conjugacies $h_{g_0,g_1}$ satisfies $h_{g_0,g_1}(0) = 0$, $h_{g_0,g_1}(1) = 1$ and $h_{g_0,g_1}(\infty) = \infty$) and the non-existence of invariant line fields supported on the filled-in Julia set of $f^\star$ imply that all these maps are hybrid conjugated with $f^\star$. But this implies that the subexponential speed of convergence given by Proposition 2.1 is impossible, since by Theorem 1 the maps in the hybrid class of $f$ converges to $f^\star$ exponentially fast.

**Case ii.** There exists a connected complex analytic curve of fixed points which contains $f$. We will apply essentially the same argument used in Case i: Note that in a similar way we can prove that all these fixed points of the operator $\mathcal{R}^2$ are polynomial-like maps which are infinitely renormalizable: in particular their filled-in Julia sets have empty interior and all their periodic points are repelling. So there are not bifurcations of periodic points in this curve of fixed points. Use the $\lambda$-lemma [MSS] to conclude that all these fixed points are quasiconformally conjugated (the argument is as in Case i). Since the fixed point $f^\star$ does not support invariant line fields in its filled-in Julia set, we conclude that all these fixed points are hybrid conjugated, which is impossible, since iterations of maps in the hybrid class of $f^\star$ converges to the fixed point $f^\star$.

So we concluded that $f^\star$ must be a hyperbolic fixed point with codimension one stable manifold.

## Appendix: Outline of Hakim's proof

To convince the reader of the existence of parabolic petals for semi-attractive compact operators in Banach spaces, we will give an outline of Hakim's proof of the existence of parabolic domains: we do not claim any sort of originality for ourselves in the following exposition and we refer to the quite clear work [H] for details. We will use the notation introduced in Section 2.1. Consider a complex analytic operator $T$ with a semi-attractive fixed point $0$. Assume $DT_0 \cdot v = v$, $v \neq 0$. In the following lines, we will identify $B$ with $\mathbb{C} \times E^s$ by the isomorphism $(x, y) \to x \cdot v + y$.

By the Stable Manifold Theorem for compact operators (see Mañé [M]), for $\delta > 0$ and $\epsilon > 0$ small the set
$$W^s_{\delta,\epsilon} = \{x \colon \exists\, C \text{ s.t. } |T^n x| < \delta \text{ and } |T^n x| \leq C(1-\epsilon)^n, \text{ for } n \geq 0\}$$
is a codimension one complex analytic manifold. More precisely, there exists a holomorphic function $\psi \colon V \to \mathbb{C}$, where $V$ is a neighborhood of $0$ on $E^s$, with $D\psi(0) = 0$, so that
$$W^s_{\delta,\epsilon} = \{(\psi(y), y) \colon y \in V\}.$$
In particular, after the local biholomorphic changes of variables

(3)
$$\begin{aligned} X &= x + \psi(y) \\ Y &= y \end{aligned}$$

it is possible to represent $T$ as $T \colon \mathbb{C} \times E^s \to \mathbb{C} \times E^s$, where $T(x, y) = (x', y')$, with

(4)
$$\begin{aligned} x' &= F(x, y) = a_1(y)x + O_y(x^2) \\ y' &= G(y) + xh(x, y) \end{aligned}$$

where $G$ is a (compact) contraction around $0$ and $a_1(0) = 1$. After the local biholomorphic change of variables

(5)
$$\begin{aligned} X &= v(y)x \\ Y &= y \end{aligned}$$

where
$$v(y) := \Pi_{i \geq 0}\ a_1(G^i(y)),$$
we can assume that $a_1 \equiv 1$.

Note that, for every $n$, $T$ has the form

(6)
$$\begin{aligned} x' &= F(x, y) = x + \sum_{2 \leq i \leq n} a_i(y)x^n + O_y(x^{n+1}) \\ y' &= G(y) + xh(x, y) \end{aligned}$$

where $G$ is a (compact) contraction around $0$. We claim that we can assume, after certain biholomorphic changes of variables, that $a_2, a_3, \cdots, a_n$ do not depend on $y$. Indeed, assume by induction that $T$ can put in the form

(7)
$$\begin{aligned} x' &= F(x, y) = x + \sum_{2 \leq i \leq n} \tilde{a}_i x^n + \tilde{a}_{n+1}(y)x^{n+1} + O_y(x^{n+2}) \\ y' &= G(y) + xh(x, y) \end{aligned}$$

Then after the local change of variables

(8)
$$X = x + v(y)x^{n+1}$$
$$Y = y$$

where $v(y) := \sum_{i \geq 0}(\tilde{a}_{n+1}(G^i(y)) - \tilde{a}_{n+1}(0))$, $T$ will have the form

(9)
$$x' = F(x,y) = x + \sum_{2 \leq i \leq n} \tilde{a}_i x^n + \tilde{a}_{n+1}(0)x^{n+1} + \tilde{a}_{n+2}(y)x^{n+2} + O_y(x^{n+2})$$
$$y' = G(y) + xh(x,y).$$

Now we are going to introduce the concept of multiplicity of the fixed point 0 for transformations on the form of Eq. (6). By the implicit function theorem, for each transformation in that form there exists a complex analytic curve $y \colon U \subset \mathbb{C} \to E^s$, with $0 \in U$, which is the unique solution for the equation

$$y(x) = G(y(x)) + xh(x, y(x)).$$

Consider the function $q \colon U \to \mathbb{C}$ defined by

$$q(x) := F(x, y(x)) - x.$$

The **multiplicity** of $T$ at 0 is defined as the order of $q$ at 0. Note that the multiplicity of $T$ at 0 is finite if and only if 0 is an isolated fixed point and infinity if and only if $q(x)$ vanishes everywhere and $(x, y(x))$ is a complex analytic curve of fixed points for $T$ (which contains all the fixed points in a neighborhood of 0). Moreover, if $T$ has the form Eq. (7), with $\tilde{a}_2 = \cdots = \tilde{a}_{n-1} = 0$ and $\tilde{a}_n \neq 0$, then the multiplicity of $T$ is exactly $n$.

Consider a transformation $T$ as in Eq. (6) and biholomorphic change of variables $W(x, y) = (X, Y)$ of the type

(10)
$$X = x + v(y)x^k$$
$$Y = y$$

where $v$ is a holomorphic function and $k > 1$. Then $W^{-1} \circ T \circ W$ has also the form in Eq. (6). Moreover

**Proposition 4.1.** *The multiplicity of $W^{-1} \circ T \circ W$ at 0 is equal to the multiplicity of $T$ at 0.*

*Proof.* (suggested by M. Lyubich) Assume that it is finite (otherwise the invariance is trivial): then 0 is an isolated fixed point. Consider the one-parameter family of change of variables $W_\lambda$ defined by

(11)
$$X = x + \lambda v(y)x^k$$
$$Y = y$$

Then $W_\lambda^{-1} \circ T \circ W_\lambda$ has the form

(12)
$$x' = F_\lambda(x, y) = x + O_{y,\lambda}(x^2)$$
$$y' = G_\lambda(y) + xh_\lambda(x, y)$$

Note that we can choose $\delta_0$ small enough so that for all $|\lambda| \leq 1$, 0 is the unique fixed point for $W_\lambda^{-1} \circ T \circ W_\lambda$ on $\{(x,y), |x| \leq \delta_0, |y| \leq \delta_0\}$. Moreover, by the implicit function theorem and the compactness of $\{\lambda\colon |\lambda| \leq 1\}$ there exists a holomorphic function $y_\lambda(x) = y(\lambda, x)$, defined on

$$\{\lambda\colon |\lambda| < 1 + \delta_1\} \times \{x\colon |x| < \delta_2\}$$

so that

$$y_\lambda(x) = G_\lambda(y) + xh_\lambda(x, y_\lambda(x))$$

Choosing $\delta_1, \delta_2$ small enough, for each $\lambda$ the point $0$ is the unique solution for the equation

$$q_\lambda(x) := F_\lambda(x, y_\lambda(x)) - x = 0$$

on $\{x\colon |x| \leq \delta_1\}$. By Rouché's Theorem, $ord_0\ q_\lambda$ does not depends on $\lambda$. □

Assume that $T$ has the form of Eq. (6) and finite multiplicity $n$. After appropriated changes of variables, we can assume that $a_2, \ldots, a_{2n-1}$ does not depends on $y$. Since the multiplicity is invariant by the above changes of variables, we conclude that $a_2 = \cdots = a_{n-1} = 0$ and $a_n \neq 0$. Doing appropriated changes of variables in the form of Eq. (10) (indeed, in this case $v$ does not depends on $y$) and replacing the coordenate $x$ by $\theta x$, for some $\theta \neq 0$, if necessary, it is possible to put $T$ in the form

$$
\begin{aligned}
(13)\quad & x' = x - \frac{1}{n-1}x^{n-1} + ax^{2(n-1)} + O_y(|x|^{2(n-1)+1}) \\
& y' = G(y) + xh(x, y).
\end{aligned}
$$

Under the above form, the set

$$P_{R,\rho} = \{(x, y)\colon |x^{n-1} - \frac{1}{2R}| < \frac{1}{2R}\ and\ |y| < \rho\}$$

is a parabolic domain, provided $R$ and $\rho$ are small enough: here Hakim's proof is very similar to the one-dimensional situation: make the "change of variables"

$$
(14)\quad \begin{aligned} X &= x^{n-1} \\ Y &= y \end{aligned}
$$

and

$$
(15)\quad \begin{aligned} X &= 1/x \\ Y &= y \end{aligned}
$$

to put $T$ in the form

$$
(16)\quad \begin{aligned}
x' &= x + 1 + c\frac{1}{x} + O_y(\frac{1}{|x|^{1+1/(n-1)}}) \\
y &= G(y) + O_y(\frac{1}{|x|^{1/(n-1)}}).
\end{aligned}
$$

and now the proof is easy.

Department of Mathematics, University of Toronto, Toronto, Ontario, Canada.

*Home page:* **www.math.sunysb.edu/˜smania/**
*E-mail address*: `dsbrandao@yahoo.com`